\documentclass{siamart0516}

\usepackage{bbm}
\newtheorem{remark}[theorem]{Remark}

\usepackage{mathtools}

\usepackage{lipsum}
\usepackage{amsfonts}
\usepackage{graphicx}
\usepackage{epstopdf}
\usepackage{algorithmic}
\ifpdf
  \DeclareGraphicsExtensions{.eps,.pdf,.png,.jpg}
\else
  \DeclareGraphicsExtensions{.eps}
\fi

\newcommand{\TheTitle}{Analysis and optimization of weighted ensemble sampling} 
\newcommand{\TheAuthors}{D. Aristoff}

\headers{\TheTitle}{\TheAuthors}

\title{{Analysis and optimization of weighted ensemble sampling}\thanks{This work was supported by the
National Science Foundation 
via the award
NSF-DMS-1522398.}}

\author{
  David Aristoff\thanks{Colorado State University}
    (\email{aristoff@math.colostate.edu})}

\usepackage{amsopn}

\begin{document}

\maketitle

\begin{abstract}
We give a mathematical 
framework for weighted ensemble (WE)
sampling, a binning and 
resampling technique for 
efficiently computing 
probabilities in molecular 
dynamics. We prove that WE sampling 
is unbiased in a very general setting 
that includes adaptive binning. 
We show that when WE is 
used for stationary calculations 
in tandem with a coarse model, the 
coarse model can 
be used to optimize the 
allocation of replicas 
in the bins. 
\end{abstract}



\section{Introduction}

This article concerns a 
resampling procedure, 
called weighted ensemble 
(WE), 
for Markov chains. WE consists 
of simulating some replicas of a Markov 
chain $(X_p)_{p \ge 0}$ 
and resampling 
from the replicas at 
certain time intervals. 
In the literature, WE sampling~\cite{dan2,darve,huber,suarez2,suarez,dan1} 
usually refers to a resampling technique  designed so that 
the replicas of $(X_p)_{p \ge 0}$ are 
evenly 
distributed throughout state space.  
This is usually achieved by dividing 
state space into bins  
and resampling in each bin so that 
the number of replicas
therein remains roughly constant.
The replicas carry probabilistic 
weights so 
that the resulting statistical 
distribution is unbiased. 
This distribution can be used, in 
principle, to estimate 
any function of 
$(X_p)_{p \ge 0}$ at a fixed time~\cite{dan1}. In practice, 
the quality of such estimates 
depends on the choice of bins and 
number of replicas maintained 
in each bin, among 
other factors. 
Below, 
we will usually refer to a replica 
as a {\em particle} and 
to resampling as {\em selection}, 
following convention in the 
mathematical literature.

Since WE is simply a resampling 
technique, it can be understood 
in a framework similar to 
that of particle filters or 
sequential Monte Carlo (SMC). 
For a review 
of standard SMC
methods, see for instance the 
textbook~\cite{del_moral}, the articles~\cite{del_moral2,del_moral3} or the 
compilation~\cite{doucet}. (See also~\cite{cerou,cerou2} 
for a related method.) 
We emphasize that WE does {not} 
fall into the SMC framework of~\cite{del_moral}, 
as there are no underlying potential 
functions or Gibbs-Boltzmann measures 
defining the selection step. We consider a very general 
framework for WE in which, contrary to the SMC/Feynman-Kac formalism (see~\cite{del_moral}), the rule for 
killing or splitting replicas is essentially arbitrary. This 
means that WE requires an independent analysis. 

The main contributions of this 
article are as follows. First, we 
give a definition of WE sampling that 
is bin-free and generalizes descriptions currently 
in the literature (Section~\ref{sec:notation}). We prove 
WE is unbiased in this setting, 
which allows for adaptive selection 
procedures~\cite{dan1} 
(Section~\ref{sec:martingale}). 
Then, we give simple 
formulas for the variance of WE 
and show how, in principle, 
the variance can be minimized 
under a constraint on the number 
of particles (Sections~\ref{sec:martingale}-\ref{sec:min}). 
In practice, the variance formula contains 
terms that may not be efficiently 
computable.
However, we show how a coarse model 
can be used to guide WE sampling 
to minimize variance in computations 
of fixed time as well 
as stationary averages of $(X_p)_{p \ge 0}$ 
(Sections~\ref{sec:binning}-\ref{sec:example}).

Our interest in WE arises 
from longstanding problems in computational 
chemistry. In 
this setting, $(X_p)_{p \ge 0}$ is obtained from a 
discretization of some 
stochastic 
molecular dynamics (MD). 
MD simulations  
have proven useful for understanding 
many chemical and biological 
processes; see~\cite{schlick} 
for an overview. However, 
such simulations are limited 
by time scale 
separation. Many phenomena of 
interest occur at the laboratory 
time scale of microseconds, 
while MD simulations have time 
steps that correspond to femtoseconds. 
In this case, straightforward MD simulations 
are not practical. Many 
methods exist for extending 
the time scale of MD simulations; 
we do not attempt to give a 
review of them here. WE is 
one of several methods for 
extending the time scale of 
simulations in models with 
rough energy landscapes. Methods 
that are related in 
scope and design include Exact Milestoning~\cite{elber,faradjian}, 
Non-Equilibrium Umbrella Sampling~\cite{warmflash,vankoten}, Trajectory Tilting~\cite{vanden-eijnden}, Transition Interface Sampling~\cite{van-erp}, 
Forward Flux Sampling~\cite{allen}, 
and Boxed Molecular Dynamics~\cite{glowacki}. 
See for instance~\cite{aristoff,darve_ryu} for 
review and comparison 
of these methods. We will comment 
on Exact Milestoning in the Appendix below. 

While WE 
can be used with a broad 
range of stochastic processes, 
when the process is 
time homogenous and 
Markovian -- as 
in many models of MD, such 
as Langevin dynamics --
WE may be used to efficiently compute 
dynamical quantities like 
reaction rates 
using a long time or 
stationary average~\cite{dan2,suarez}. 
These computations rely on 
Hill's relation~\cite{hill}, 
which we generalize in the Appendix below. 
From Hill's relation, 
obtaining reaction 
rates  
requires a calculation using
the stationary distribution 
of a nonreversible process. 
 
To speed up the stationary 
calculation, WE is combined with a 
preconditioning step~\cite{dan2,darve} 
in which a Markov state model (MSM)~\cite{sarich,schutte} is used 
to approximate the stationary 
distribution. This is 
sometimes called accelerated 
WE~\cite{darve}. 
Accelerated WE begins with 
particles evenly distributed 
in space, with weights chosen 
to match the stationary distribution 
of the MSM. The particles are 
then allowed to relax according 
to their exact dynamics, with WE 
sampling 
ensuring that the particles 
remain evenly distributed in 
state  
space. 
We show below that information 
from the MSM can be used to 
optimize the WE sampling in 
this relaxation step, in the 
sense that the variance in 
the appropriate stationary 
calculation is minimized. 
This optimization requires 
an adaptive number of 
particles per bin, in 
contrast with traditional 
WE sampling. We show 
in a simple model 
that this adaptive 
sampling can be significantly 
better than traditional 
WE sampling.

This article is organized as 
follows. In Section~\ref{sec:notation}, 
we define a WE process in a general setting and give 
an algorithm for 
WE sampling.
In Section~\ref{sec:martingale}, 
we introduce a martingale framework for WE sampling in this setting. 
We use the framework to prove the sampling 
is unbiased and obtain formulas for the variance. In Section~\ref{sec:min}, 
we show how to minimize the variance 
under a constraint on the total number 
of particles.
In Section~\ref{sec:binning} we consider WE sampling based on  
binning techniques. 
In Section~\ref{sec:adapt} we show how 
adapting the binning 
to a coarse 
model for $(X_p)_{p\ge 0}$ can be 
used to minimize variance, and 
in Section~\ref{sec:stationary} we 
apply these ideas to computing 
stationary averages.
In Section~\ref{sec:example}, 
we use a simple model 
to compare adaptive WE 
to traditional WE and 
naive sampling. In the Appendix, 
we prove a generalization of 
the Hill relation and discuss 
connections to Exact Milestoning.

\section{Notation and assumptions}\label{sec:notation}

Throughout, $(X_p)_{p \ge 0}$ is a time homogeneous 
Markov chain with values in a 
measurable space $(E,{\mathcal E})$ 
and transition kernel $K$. We write $\sim$ 
to denote equality in law of random 
variables or processes, and 
${\mathbb E}$ and ${\mathbb P}$ 
for various expectations and probabilities. 
When $\zeta$ is a probability measure 
on $(E,{\mathcal E})$, superscripts such as ${\mathbb E}^\zeta$ 
or ${\mathbb P}^\zeta$ represent 
processes with initial distribution $\zeta$, with ${\mathbb E}^x$ 
or ${\mathbb P}^x$ indicating the processes start at the 
point $x$.
Sets and functions will be assumed 
measurable without explicit mention. 
For a measure $\zeta$ on $(E,{\mathcal E})$ 
and bounded $f:E \to {\mathbb R}$, we write $\zeta(f) = \int f\,d\zeta$ for 
the integral of $f$ with respect 
to $\zeta$. We also write $\zeta K(dy) = \int K(x,dy)\zeta(dx)$ for left action of $K$, 
and $Kf(x) = \int K(x,dy)f(y)$ for 
the right action. In particular, 
\begin{equation}\label{eq_LR_action}
\zeta K^n f = {\mathbb E}^\zeta[f(X_n)]
\end{equation}
where throughout $K$ is the 
Markov transition kernel of $(X_p)_{p \ge 0}$.
If $S$ is a set, we write $\#S$ for the number 
of elements of $S$.

We study a certain class of 
sequential Monte Carlo (SMC)
methods for sampling 
$(X_p)_{p \ge 0}$ described in Definition~\ref{def_PS} 
below.
We begin with an informal description 
of the procedure.
Consider a process 
consisting of particles 
in $E$ and weights in ${\mathbb R}^+ = [0,\infty)$. 
The initial particles 
all have the same distribution as $X_0$.
At each time $p$, some of the particles 
are selected, or copied, and others 
are thrown away, or killed. The 
selected particles then 
mutate according to the evolution law of $(X_p)_{p \ge 0}$. 
(We often refer to selected 
particles as {\em children} 
and the particles from which 
they were copied as {\em parents}.)
The selected points and 
weights are chosen to yield unbiased 
estimators for the law of $(X_p)_{p \ge 0}$. 
This is ensured by setting a child's 
weight to be its parent's 
weight divided by the expected 
number of times the parent is selected.
Writing $\xi_p^j$ and $\omega_p^j$ 
for the particles and weights at time $p$, 
and using the symbol ``$\,\,\hat{}\,\,$'' to indicate selected 
particles and weights,
we make this precise as follows.
\vskip5pt
\begin{definition}\label{def_PS}
A {\em weighted ensemble} (WE) consists of particles and weights 
\begin{equation*}
(\xi_p^j,\omega_p^j)_{p \ge 0}^{j=1,\ldots,N_p}, \qquad ({\hat \xi}_p^i,{\hat\omega}_p^i)_{p\ge 0}^{i=1,\ldots,{\hat N}_p}
\end{equation*}
with 
values in $\cup_{n=1}^\infty (E \times {\mathbb R}^+)^n$, 
selection rules $$(C_p^j)_{p \ge 0}^{j=1,\ldots,N_p}$$ with values in $\cup_{n=1}^{\infty}({\mathbb N}\cup \{0\})^n$, and associated filtrations 
\begin{align*}
{\mathcal F}_p &= \sigma\left((\xi_q^j,\omega_q^j)_{0 \le q \le p}^{j=1,\ldots,N_q},(C_q^j)_{0 \le q \le p-1}^{j=1,\ldots,N_q},({\hat \xi}_q^i,{\hat \omega}_q^i)_{0 \le q \le p-1}^{i=1,\ldots,{\hat N}_q}\right) \\
{\hat {\mathcal F}}_p&= \sigma\left((\xi_q^j,\omega_q^j)_{0 \le q \le p}^{j=1,\ldots,N_q},(C_q^j)_{0 \le q \le p}^{j=1,\ldots,N_q},({\hat \xi}_q^i,{\hat \omega}_q^i)_{0 \le q \le p}^{i=1,\ldots,{\hat N}_q}\right)
\end{align*}
which together satisfy (A1)-(A4) below for each $p \ge 0$.
\vskip5pt
\begin{itemize}
\item[(A1)] $N_0>0$ is constant, and for 
$j=1,\ldots,N_0$, 
$\xi_0^j\sim X_0$, $\omega_0^j = 1/N_0$.
\vskip2pt
\item[(A2)] Each child
${\hat \xi}_p^i$ is associated 
to a parent $\xi_p^{\alpha(i)}$. 
With 
\begin{equation*}
C_p^j = \#\{i: \alpha(i) = j\}
\end{equation*}
the number of children of $\xi_p^j$, we have ${\mathbb E}[C_p^j|{\mathcal F}_p] > 0$,
\begin{equation*}
{\hat N}_p = \sum_{j=1}^{N_p} C_p^j,
\qquad \text{and} \qquad
{\hat \omega}_p^i = \frac{\omega_p^j}{{\mathbb E}\left[C_p^j|{\mathcal F}_p\right]}, \quad \text{ if }\alpha(i) = j.
\end{equation*}
\vskip2pt
\item[(A2')] Conditionally on ${\mathcal F}_p$,  $C_p^1,\ldots,C_p^{N_p}$ are 
uncorrelated.
\vskip2pt
 \item[(A3)] 
 $N_{p+1} = {\hat N}_p$ and 
 $\omega_{p+1}^i = {\hat \omega}_p^i$ 
 for $i = 1,\ldots,{\hat N}_p$.
 \vskip2pt
 \item[(A4)] Conditionally on ${\hat {\mathcal F}}_p$,
 $\xi_{p+1}^1,\ldots, \xi_{p+1}^{N_{p+1}}$ are independent with 
 \begin{equation*}
 {\mathbb P}[\xi_{p+1}^i \in dx] = K({\hat \xi}_p^i,dx).
 \end{equation*}
\end{itemize}

\end{definition}

It is convenient to view a WE through the 
following diagram:
\begin{align*}
&\{\xi_p^j\}^{j=1,\ldots,N_p}
\xrightarrow{\textup{selection}} 
\{{\hat \xi}_p^i\}^{i=1,\ldots,{\hat N}_p}
\xrightarrow{\textup{mutation}} 
\{{\xi}_{p+1}^j\}^{j=1,\ldots,N_{p+1}},\\
&\{\omega_p^j\}^{j=1,\ldots,N_p}
\xrightarrow{\textup{selection}} 
\{{\hat \omega}_p^i\}^{i=1,\ldots,{\hat N}_p}
\xrightarrow{\textup{mutation}} 
\{{\omega}_{p+1}^j\}^{j=1,\ldots,N_{p+1}}.
\end{align*}
The filtration ${\mathcal F}_p$ (resp. ${\hat{\mathcal F}}_p$) represents the information from the WE process up to time $p$, not including the selection step 
(resp. up to time $p$, including the 
selection step).  
We write $\alpha(i)$ for 
the index of the parent 
particle of the $i$th selected particle. Thus, 
\begin{equation*}
\alpha(i) = j \Longrightarrow {\hat \xi}_p^i = \xi_p^j.
\end{equation*}
(Of course $\alpha$ depends on $p$, but 
we do not make this explicit.) Also,
\begin{equation*}
C_p^j = \#\{i:\alpha(i) = j\} = \text{number of times }\xi_p^j \text{ is selected}.
\end{equation*}
The $C_p^j$, $j=1,\ldots,N_p$, can depend on the entire history of 
the process. 
We assume in (A2') that they are uncorrelated conditionally on 
the past so that we can 
obtain a relatively simple explicit formula 
for variance in Theorem~\ref{thm_martingale} 
below. This assumption is 
only needed 
for the variance. 
Indeed, the proof of Theorem~\ref{thm_martingale} below shows that 
(A2') is not required for unbiased WE sampling; see the remarks after the proof of
Theorem~\ref{thm_martingale}.

Note that, by 
(A2), the weight of a selected 
particle 
is simply the weight of its parent 
particle 
divided by the expected number 
of times the parent is selected. 
We assume the expected 
number of times a parent is selected 
is positive, so that each particle 
has a positive probability to survive.

Assumption
(A1) says that the initial collection 
of particles and weights is 
chosen according to the distribution of $X_0$. 
Notice we do not 
require that the $\xi_0^i$'s are 
independent, so they can 
be generated by, for example, 
Markov chain Monte Carlo or 
other sequential samplers.

\begin{algorithm}
\caption{A WE sampler}\label{alg0} 
Choose initial points and weights 
$(\xi_0^j,\omega_0^j)^{j=1,\ldots,N_0}$ 
according to the distribution of $X_0$ in the sense of (A1). Then iterate over $p \ge 0$ until time $\tau_{kill} := \inf\{p\ge 0:N_p = 0\}$: 
\vskip5pt
\begin{itemize}
\item[1.] For $j=1,\ldots,N_p$, choose a 
number $C_p^j$ of times to select 
particle $\xi_p^j$. Let ${\hat \xi}_p^i$, 
$i=1,\ldots,{\hat N}_p$ be the collection 
of selected particles, with 
${\hat N}_p = \sum_{j=1}^{N_p} C_p^j$.
\item[2.] Assign the weight ${\hat \omega}_p^i = \frac{\omega_p^j}{{\mathbb E}[C_p^j|{\mathcal F}_p]}$ to ${\hat \xi}_p^i$, if $\alpha(i) = j$.
\item[3.] Set $N_{p+1} = {\hat N}_p$ and 
$\omega_{p+1} = {\hat \omega}_p^i$ 
for $i=1,\ldots,{\hat N}_p$.
\item[4.] Evolve the particles ${\hat \xi}_p^i$, $i=1,\ldots,{\hat N}_p$, independently
according to the law of $(X_p)_{p\ge 0}$ to get the next generation $\xi_{p+1}^j$,
$j=1,\ldots,N_{p+1}$ of particles.
\end{itemize}
\vskip5pt
Steps 1-2 correspond to selection, 
and 3-4 to evolution. In the above, the $C_p^j$'s are usually independent of each other, given the current state of the algorithm,  
and they can depend on the entire history of the algorithm 
up to time $p$. In Step 2, 
${\mathbb E}[C_p^j|{\mathcal F}_p]$ 
represents the expected 
value of $C_p^j$ given that history. 
We assume ${\mathbb E}[C_p^j|{\mathcal F}_p]>0$, that is, every particle has a positive survival 
probability.
\end{algorithm}

Assumption (A3) says that the 
weights defined in the selection 
step will be assigned to the 
particles in the next generation.
The condition 
(A4) states that the next 
generation of particles 
mutates from the selected 
particles using the evolution 
law of $(X_p)_{p \ge 0}$, where these particles 
evolve 
independently from each other. 

For clarity, we give an algorithm 
for sampling a WE process; see Algorithm~\ref{alg0}.

We will show in Theorem~\ref{thm_martingale} 
below that a WE in the sense of Definition~\ref{def_PS} 
is an unbiased estimator 
for the law of $(X_p)_{p \ge 0}$. To make this precise 
we introduce the following notation.
At time $p$, a WE defines 
empirical distributions
\begin{equation}\label{eq_emp_distr}
\eta_p = \sum_{j=1}^{N_p} \omega_p^j 
\delta_{\xi_p^j}, 
\qquad {\hat \eta}_p = \sum_{i=1}^{{\hat N}_p} {\hat \omega}_p^i 
{\delta}_{{\hat \xi}_p^i}.
\end{equation}
These definitions make sense only up until the first time all the particles have been killed, $
\tau_{kill} = \inf\{p \ge 0:{N}_p = 0\} = \inf\{p \ge 0: {\hat N}_{p-1} = 0\}$.
We adopt the convention that 
$\eta_p \equiv 0$ and ${\hat \eta}_{p-1} \equiv 0$ if 
$p\ge \tau_{kill}$. It is important to note that $\eta_p(1) \ne 1$ in general; that is, the total weight is not 
conserved exactly. 
\vskip5pt
\begin{remark}\label{rmk1}
Often it is desirable to fix 
the average total number of particles, or 
simply the total number of particles. 
Below we consider mostly the former case, 
but here we comment briefly on the latter.

A simple population control step 
can be added to guarantee $N_p \equiv N$ 
for each $p$,
with $N$ fixed, as follows. First, assume the population control has been applied 
up to time $p$, so that $N_p = N$. Suppose 
furthermore that the selection is 
done so that
${\mathbb E}[N_{p+1}|\mathcal{F}_p] = N$.
Then 
\begin{align*}
N =  {\mathbb E}[{\hat N}_p|\mathcal{F}_p] 
= {\mathbb E}\left[\left.\sum_{j=1}^{N_p}C_p^j\right|\mathcal{F}_p\right] 
= \sum_{j=1}^N {\mathbb E}[C_p^j|{\mathcal F}_p].
\end{align*}
It follows that ${\mathbb E}[C_p^j|{\mathcal F}_p] \ge 1$ for some $j$. For this $j$ 
we may assume $C_p^j \ge 1$ with probability $1$, conditional on ${\mathcal{F}}_p$. 
Thus, we can 
assume there is no extinction 
in the selection step. 
Then, after 
the selection step, 
we can kill or copy particles 
uniformly at random to 
enforce $N_{p+1} = N$, and adjust weights accordingly. Note that 
this extra step would introduce 
correlations between the 
number of children of 
each particle, which 
would violate (A2'). So that 
we can obtain simple 
variance formulas, below we will focus on the 
case of uncorrelated 
$C_p^j$'s.
\end{remark}

\section{Martingale framework and variance}\label{sec:martingale}

Recall that $K$ is the transition kernel of $(X_p)_{p \ge 0}$, and recall the definitions of $\eta_p$ and 
${\hat \eta}_p$ from~\eqref{eq_emp_distr}. 
In this section and below, $n\ge 0$ and a bounded 
function $f:E \to {\mathbb R}$ are fixed.
For $0 \le p \le n$ define
\begin{equation*}
M_p = \eta_p K^{n-p} f, \qquad {\hat M}_p = {\hat \eta}_p K^{n-p}f,
\end{equation*} 
where by convention $K^0f = f$. Since $f$ 
is bounded, both $(M_p)_{0 \le p \le n}$ and $({\hat M}_p)_{0 \le p\le n}$ are integrable 
and square integrable.
Intuitively, $M_p$ represents starting at the distribution $\eta_p$, evolving 
forward $n-p$ time steps using $K$, 
and then evaluating $f$; see~\eqref{eq_LR_action}. Our analysis 
below is based on the following 
result. 
\vskip5pt
\begin{theorem}\label{thm_martingale}
Let assumptions (A1), (A2), (A3), and (A4) hold. Then $(M_p)_{0 \le p \le n}$ is a ${\mathcal F}_p$-martingale and 
\begin{equation}\label{eq_l2error}
{\mathbb E}[M_n^2] = {\mathbb E}[M_0^2] 
+ {\mathbb E}\left[\sum_{p=0}^{n-1} \left(
{\mathbb E}\left[(M_{p+1}-{\hat M}_p)^2|\hat{\mathcal F}_{p}\right]\right] + {\mathbb E}\left[({\hat M}_p - M_p)^2|{\mathcal F}_p\right]\right).
\end{equation}
If in addition (A2') holds, 
then with $g_p = K^{n-p}f$,
\begin{align}\begin{split}\label{eq_l2explicit}
&{\mathbb E}\left[\left.(M_{p+1}-{\hat M}_{p})^2\right|{\hat {\mathcal F}}_{p}\right] = \sum_{i=1}^{{\hat N}_p} ({\hat \omega}_p^i)^2[Kg_{p+1}^2({\hat \xi}_p^i) - g_p({\hat \xi}_p^i)^2] \\
&{\mathbb E}\left[\left.({\hat M}_{p}-M_{p})^2\right|{\mathcal F}_{p}\right] = \sum_{j=1}^{N_p}  ({\omega}_p^j)^2\left[\frac{{\mathbb E}[(C_p^j)^2|{\mathcal F}_p]}{{\mathbb E}[C_p^j|{\mathcal F}_p]^2}  -1\right]g_p(\xi_p^j)^2
\end{split}
\end{align}
for $0 \le p \le n-1$.
\end{theorem}
\vskip5pt
\begin{proof}
Consider a WE as in Definition~\ref{def_PS}. From (A4),
\begin{equation}\label{first}
{\mathbb E}\left[\left.
g_{p+1}(\xi_{p+1}^i)\right|
\hat{{\mathcal F}_p}\right] 
= Kg_{p+1}({\hat \xi}_p^i) = g_p({\hat \xi}_p^i).
\end{equation}
If in addition (A2') holds, then
\begin{equation}\label{second}
{\mathbb E}\left[\left.
g_{p+1}(\xi_{p+1}^i)g_{p+1}(\xi_{p+1}^k)\right|
\hat{{\mathcal F}_p}\right] =  
\begin{cases} g_p({\hat \xi}_p^i)
g_p({\hat \xi}_p^k), & i\ne k\\
K g_{p+1}^2({\hat \xi}_p^i), 
& i = k\end{cases}.
\end{equation} 
Suppose (A1), (A2), (A3) and (A4) hold. 
We show that then $$M_0,{\hat M}_0,M_1,{\hat M}_1,\ldots,{\hat M}_{n-1},M_n$$ is a martingale with respect 
to the filtration $${\mathcal F}_0 \subseteq \hat{\mathcal F}_0 \subseteq {\mathcal F}_1 \subseteq \hat{\mathcal F}_1 \subseteq \ldots \subseteq \hat{\mathcal F}_{n-1} \subseteq {\mathcal F}_n.$$ That is, we show 
that for $0 \le p \le n-1$,
\begin{equation*}
{\mathbb E}[M_{p+1}|\hat{\mathcal F}_p] = \hat{M}_p, 
\qquad {\mathbb E}[{\hat M}_p|{\mathcal F}_p] = M_p.
\end{equation*}
The fact that $M_p$ is a martingale then follows from 
\begin{equation*}
{\mathbb E}[M_{p+1}|{\mathcal F}_p] = 
{\mathbb E}[{\mathbb E}[M_{p+1}|\hat{\mathcal F}_p]|{\mathcal F}_p] = {\mathbb E}[{\hat M}_p|{\mathcal F}_p] 
= M_p,
\end{equation*}
and equation~\eqref{eq_l2error} follows from the Doob decomposition. Since $g_p = K^{n-p}f$, 
\begin{equation*}
M_p = \sum_{j=1}^{N_p} \omega_p^j g_p(\xi_p^j), \qquad 
{\hat M}_p = \sum_{i=1}^{{\hat N}_p} {\hat\omega}_p^i g_p({\hat \xi}_p^i).
\end{equation*}
So by (A3) and~\eqref{first},
\begin{align*}
{\mathbb E}[M_{p+1}|\hat{\mathcal F}_p] 
&= {\mathbb E}\left[\left.\sum_{i=1}^{N_{p+1}}\omega_{p+1}^i g_{p+1}(\xi_{p+1}^i)\right|\hat{\mathcal F}_p\right] \\
&= \sum_{i=1}^{{\hat N}_p} {\hat \omega}_p^i {\mathbb E}[g_{p+1}(\xi_{p+1})|\hat{\mathcal F}_p] \\
&= \sum_{i=1}^{{\hat N}_p} {\hat \omega}_p^i g_p({\hat \xi}_p^i) = {\hat M}_p.
\end{align*}
Also, by (A2),
\begin{align*}
{\mathbb E}[{\hat M}_p|{\mathcal F}_p] 
&= {\mathbb E}\left[\left.\sum_{i=1}^{{\hat N}_{p}}{\hat \omega}_{p}^i g_{p}({\hat \xi}_{p}^i)\right|{\mathcal F}_p\right] \\
&= \sum_{j=1}^{N_p}{\mathbb E}\left[\left.\sum_{i:\alpha(i) = j}{\hat \omega}_{p}^i g_{p}({\hat \xi}_{p}^i)\right|{\mathcal F}_p\right] \\
&= \sum_{j=1}^{N_p}\frac{\omega_p^j}{{\mathbb E}[C_p^j|{\mathcal F}_p]}{\mathbb E}\left[\left.\sum_{i:\alpha(i) = j} g_{p}({\hat \xi}_{p}^i)\right|{\mathcal F}_p\right] \\
&= \sum_{j=1}^{N_p} \omega_p^j g_p(\xi_p^j) = M_p.
\end{align*}
It remains to establish~\eqref{eq_l2explicit}. Suppose in 
addition (A2') holds. 
By~\eqref{second} and (A3),
\begin{align*}
{\mathbb E}[M_{p+1}^2|\hat{\mathcal F}_p] &= 
{\mathbb E}\left[\left.\sum_{i,k=1}^{N_{p+1}} \omega_{p+1}^i \omega_{p+1}^k g_{p+1}(\xi_{p+1}^i)g_{p+1}(\xi_{p+1}^k)\right|\hat{\mathcal F}_p\right] \\
&= \sum_{i,k=1}^{{\hat N}_p} {\hat \omega}_p^i{\hat \omega}_p^k {\mathbb E}\left[\left.g_{p+1}(\xi_{p+1}^i)g_{p+1}(\xi_{p+1}^k)\right|\hat{\mathcal F}_p\right]\\
&= \sum_{\substack{i,k=1\\i \ne k}}^{{\hat N}_p} {\hat \omega}_p^i {\hat \omega}_p^k g_p({\hat \xi}_p^i)g_p({\hat \xi}_p^k) + \sum_{i=1}^{{\hat N}_p} 
({\hat \omega}_p^i)^2 Kg_{p+1}^2({\hat \xi}_p^i).
\end{align*}
Subtracting ${\hat M}_p^2$ from this gives 
\begin{equation*}
{\mathbb E}[(M_{p+1}-{\hat M}_p)^2|\hat{\mathcal F}_p] 
= {\mathbb E}[M_{p+1}^2|\hat{\mathcal F}_p] - {\hat M}_p^2 = \sum_{i=1}^{{\hat N}_p} ({\hat \omega}_p^i)^2[Kg_{p+1}^2({\hat \xi}_p^i) - g_p({\hat \xi}_p^i)^2].
\end{equation*}
Next, notice that, with $\beta_p^j = {\mathbb E}[C_p^j|{\mathcal F}_p]$, by (A2),
\begin{align*}
{\mathbb E}[{\hat M}_p^2|{\mathcal F}_p] &= {\mathbb E}\left[\left. \sum_{i,k=1}^{{\hat N}_p} {\hat \omega}_p^i {\hat \omega}_p^kg_p({\hat \xi}_p^i)g_p({\hat \xi}_p^k)\right|{\mathcal F}_p\right] \\
&= \sum_{j,\ell=1}^{N_p} \frac{\omega_p^j\omega_p^\ell}{\beta_p^j\beta_p^\ell}{\mathbb E}\left[\left.\sum_{i,k:\alpha(i)=j,\alpha(k)=\ell} g_p({\hat \xi}_p^i)g_p({\hat \xi}_p^k)\right|{\mathcal F}_p\right],
\end{align*}
Summing over $j\ne \ell$ and using (A2'), we get 
\begin{equation*}
\sum_{\substack{j,\ell=1 \\j\ne \ell}}^{N_p} \frac{\omega_p^j\omega_p^\ell}{\beta_p^j\beta_p^\ell}{\mathbb E}\left[\left.\sum_{i,k:\alpha(i)=j,\alpha(k)=\ell} g_p({\hat \xi}_p^i)g_p({\hat \xi}_p^k)\right|{\mathcal F}_p\right] = 
\sum_{\substack{j,\ell=1 \\j\ne \ell}}^{N_p} 
\omega_p^j \omega_p^\ell g_p(\xi_p^j)g_p(\xi_p^\ell),
\end{equation*}
and summing over $j = \ell$, with $\gamma_p^j = {\mathbb E}[(C_p^j)^2|{\mathcal F}_p]$, we have 
\begin{equation*}
\sum_{j=1}^{N_p} \left(\frac{\omega_p^j}{\beta_p^j}\right)^2{\mathbb E}\left[\left.\sum_{i,k:\alpha(i) = j, \alpha(k) = j} g_p({\hat \xi}_p^i)g_p({\hat \xi}_p^k)\right|{\mathcal F}_p\right] = 
\sum_{j=1}^{N_p} \left(\frac{\omega_p^j}{\beta_p^j}\right)^2\gamma_p^j g_p(\xi_p^j)^2.
\end{equation*}
Combining the last three displays, 
\begin{equation*}
{\mathbb E}[{\hat M}_p^2|{\mathcal F}_p] = 
\sum_{\substack{j,\ell=1 \\j\ne \ell}}^{N_p} 
\omega_p^j \omega_p^\ell g_p(\xi_p^j)g_p(\xi_p^\ell) + 
\sum_{j=1}^{N_p} \left(\frac{\omega_p^j}{\beta_p^j}\right)^2\gamma_p^j g_p(\xi_p^j)^2.
\end{equation*}
Subtracting $M_p^2$, we get 
\begin{equation*}
{\mathbb E}[({\hat M}_p-M_p)^2|{\mathcal F}_p] 
= {\mathbb E}[{\hat M}_p^2|{\mathcal F}_p] - M_p^2 
= \sum_{j=1}^{N_p}  \left[\left(\frac{\omega_p^j}{\beta_p^j}\right)^2\gamma_p^j  - (\omega_p^j)^2\right]g_p(\xi_p^j)^2.
\end{equation*}
\end{proof}
\vskip5pt
Below, we will repeatedly refer to the functions 
\begin{equation*}
g_p = K^{n-p}f
\end{equation*}
from the proof above, so we record the definition here again for convenience.
We note that Theorem~\ref{thm_martingale} shows WE is unbiased, as follows. 
Since
$(M_p)_{0 \le p \le n}$ is a martingale, ${\mathbb E}[M_n] = {\mathbb E}[M_0]$. Moreover,
(A1) implies ${\mathbb E}[M_0] = {\mathbb E}[f(X_n)]$. This means that
\begin{equation*}
{\mathbb E}\left[\sum_{j=1}^{N_n} \omega_n^j f(\xi_n^j)\right] = {\mathbb E}\left[M_n\right] = {\mathbb E}[f(X_n)].
\end{equation*}
The proof of Theorem~\ref{thm_martingale} 
shows that this equation does not 
require assumption (A2').  
Notice also that~\eqref{eq_l2error} leads 
to a formula for the $L^2$ sampling error, or variance, via 
\begin{align}\begin{split}\label{eq_varform}
{\mathbb E}\left[\left(\sum_{j=1}^{N_n} \omega_n^j f(\xi_n^j) - {\mathbb E}[f(X_n)]\right)^2\right] &= {\mathbb E}\left[(M_n-{\mathbb E}[f(X_n)])^2\right]\\
&= {\mathbb E}[M_n^2] 
- {\mathbb E}[f(X_n)]^2.
\end{split}
\end{align}
By Theorem~\ref{thm_martingale}, the expression in~\eqref{eq_varform} consists of a term corresponding to the variance from the initial condition, namely 
$
{\mathbb E}[M_0^2] 
- {\mathbb E}[f(X_n)]^2 = \text{Var}(M_0)$,
added to another 
term corresponding to the variance 
from the evolutions and selections, namely 
\begin{equation*}
{\mathbb E}\left[\sum_{p=0}^{n-1} \left(
{\mathbb E}\left[(M_{p+1}-{\hat M}_p)^2|\hat{\mathcal F}_{p}\right]\right] + {\mathbb E}\left[({\hat M}_p - M_p)^2|{\mathcal F}_p\right]\right).
\end{equation*}
If we assume (A2'), we can get 
nice expressions for the latter variances; 
see~\eqref{eq_l2explicit}. In~\eqref{eq_l2explicit}, we can think of the first equation 
as the variance due to {\em mutation}, 
and the second equation as the variance 
from {\em selection}. Indeed we can 
understand these as ``local variances'' 
associated to particle 
evolution and selection since 
we can rewrite
\begin{equation*}
{\mathbb E}\left[\left.(M_{p+1}-{\hat M}_{p})^2\right|{\hat {\mathcal F}}_{p}\right] = 
\sum_{i=1}^{{\hat N}_p} ({\hat \omega}_p^i)^2 
\text{Var}_{K({\hat \xi}_p^i,\cdot)}(g_{p+1})
\end{equation*} 
and 
\begin{equation*}
{\mathbb E}\left[\left.({\hat M}_{p}-M_{p})^2\right|{\mathcal F}_{p}\right] = 
\sum_{j=1}^{N_p}({\omega}_p^j)^2\frac{\textup{Var}(C_p^j|{\mathcal F}_p)}{{\mathbb E}[C_p^j|{\mathcal F}_p]^2}g_p(\xi_p^j)^2.
\end{equation*}
In the following sections
we will attempt to minimize 
these terms to produce a 
near optimal sampling strategy.

\section{Minimizing variance}\label{sec:min}

The main idea in the sections that follow 
is to use information available at time $p$ 
-- that is, the ${\mathcal F}_p$-measurable random variables -- to decide 
how to make the selections. We want to minimize the variance from 
both selection and mutation, subject 
to a constraint on the average 
total number of particles. 
Instead of trying to simultaneously
minimize both variances, we will minimize 
in two steps: first, we minimize the variance from mutation, and then, subject to the constraints 
thereby imposed, we minimize the
variance from selection.

For the variance from mutation, we have 
to condition on ${\mathcal F}_p$ to get 
an expression that depends only on ${\mathcal F}_p$-measurable random variables. Thus, using~\eqref{eq_l2explicit} and (A2),
\begin{align}\begin{split}\label{eq_mut_min}
{\mathbb E}\left[\left.(M_{p+1}-{\hat M}_{p})^2\right| {\mathcal F}_{p}\right] 
&= {\mathbb E}\left[\left.{\mathbb E}\left[\left.(M_{p+1}-{\hat M}_{p})^2\right|\hat{\mathcal F}_p\right]\right| {\mathcal F}_{p}\right]\\
&= {\mathbb E}\left[\left.
\sum_{i=1}^{{\hat N}_p} ({\hat \omega}_p^i)^2[Kg_{p+1}^2({\hat \xi}_p^i) - g_p({\hat \xi}_p^i)^2] \right|{\mathcal F}_p\right] \\
&= \sum_{j=1}^{N_p}\left(\frac{\omega_p^j}{{\mathbb E}[C_p^j|{\mathcal F}_p]}\right)^2{\mathbb E}\left[\left.\sum_{i:\alpha(i)=j} [Kg_{p+1}^2({\hat \xi}_p^i)-g_p({\hat\xi}_p^i)^2]\right|{\mathcal F}_p\right] \\
&= \sum_{j=1}^{N_p} \frac{(\omega_p^j)^2}{{\mathbb E}[C_p^j|{\mathcal F}_p]}[Kg_{p+1}^2(\xi_p^j)-g_p(\xi_p^j)^2].
\end{split} 
\end{align}
Minimizing this expression is only 
interesting if we limit the total number 
of particles. In principle, we can choose $C_p^j$'s such that this variance is minimized, given a fixed target number, $N$, of total particles. More precisely, if we demand that
\begin{equation*}
\sum_{j=1}^{N_p} {\mathbb E}[C_p^j|{\mathcal F}_p] = N,
\end{equation*}
then a Lagrange multiplier calculation shows~\eqref{eq_mut_min} is minimized by $C_p^j$'s with 
\begin{equation}\label{min_cpj}
{\mathbb E}[C_p^j|{\mathcal F}_p] = \frac{N \omega_p^j \sqrt{Kg_{p+1}^2(\xi_p^j)-g_p(\xi_p^j)^2}}{\sum_{j=1}^{N_p}\omega_p^j \sqrt{Kg_{p+1}^2(\xi_p^j)-g_p(\xi_p^j)^2}}.
\end{equation}
(provided the denominator above is nonzero).
Note that from Jensen's inequality, 
\begin{equation*}
Kg_{p+1}^2(\xi) - g_p(\xi)^2 = \text{Var}_{K(\xi,\cdot)}(g_{p+1}) \ge 0
\end{equation*}
for all $\xi \in E$, and indeed this 
expression can be understood as a ``local 
variance'' associated with mutating 
a particle $\xi$. 
Recall the variance due to selection is 
\begin{equation}\label{min_sel}
{\mathbb E}\left[\left.({\hat M}_{p}-M_{p})^2\right|{\mathcal F}_{p}\right] = \sum_{j=1}^{N_p}  (\omega_p^j)^2\left[\frac{{\mathbb E}[(C_p^j)^2|{\mathcal F}_p]}{{\mathbb E}[C_p^j|{\mathcal F}_p]^2}  - 1\right]g_p(\xi_p^j)^2.
\end{equation}  
Our minimization strategy is as follows. 
First, we choose $C_p^j$'s satisfying~\eqref{min_cpj}. 
Note that this step only determines their average values 
$$\beta_p^j = {\mathbb E}[C_p^j|{\mathcal F}_p].$$ To minimize~\eqref{min_sel} over these $C_p^j$'s, 
we take 
${\mathbb E}[(C_p^j)^2|{\mathcal F}_p]$ as small 
as possible. This is done as follows. Let $\lfloor x \rfloor$ be the 
integer part of~$x$. Then conditionally on ${\mathcal F}_p$, set each $C_p^j$ to equal either $\lfloor \beta_p^j\rfloor$ or 
$\lfloor \beta_p^j\rfloor + 1$, 
with probabilities chosen so that the mean is $\beta_p^j$.  

The problem with the above strategy is that, 
in~\eqref{min_cpj}, 
the quantities $$Kg_{p+1}^2(\xi_p^j)-g_p(\xi_p^j)^2$$ 
are not easily computable. Indeed, if they were, then 
$(X_p)_{p \ge 0}$ would be simple enough that WE is not 
needed. 
We have found that nonetheless a version of the 
above strategy 
can be made useful if we obtain coarse approximations 
for these quantities. The basic idea is to 
construct a coarse model, for instance 
a Markov State Model~\cite{sarich,schutte}, for $(X_p)_{p \ge 0}$ 
from which the $Kg_{p+1}^2(\xi_p^j)-g_p(\xi_p^j)^2$ 
can be approximated. The coarse model 
will have states that correspond to {\em bins} 
that partition $E$, and the WE process 
will be adapted to the same bins, in the sense 
that the resampling rules are tailored 
to the bin structure. We pursue these 
ideas in the following sections.

\vskip5pt
\begin{remark}
It is interesting to consider the 
limits where the time or the number of particles 
become infinite. We briefly comment on 
the latter. If we substitute the minimizing 
equation~\eqref{min_cpj} into~\eqref{eq_mut_min}, and set $$F(\xi) := \sqrt{Kg_{p+1}^2(\xi)-g_p(\xi)^2},$$ then we get
\begin{equation}\label{norm_asym}
N{\mathbb E}\left[\left.(M_{p+1}-{\hat M}_{p})^2\right| {\mathcal F}_{p}\right] 
= \left(\sum_{j=1}^{N_p} 
\omega_p^j F(\xi_p^j)\right)^2 = (\eta_p(F))^2.
\end{equation}
Intuitively, under appropriate conditions on $F$, the following particle approximation result is suggested 
by a version of the
law of large numbers for 
sufficiently weakly dependent 
random variables:
$$\eta_p(F)  \xrightarrow{a.s.} {\mathbb E}[F(X_p)]\quad \text{ as }N \to \infty.$$ 
We leave this question for future 
work; see however Section 7.4 of~\cite{del_moral} for analogous 
results in the SMC/Feyman-Kac framework. 
Taking this result for granted 
and using~\eqref{norm_asym}, we expect
${\mathbb E}[F(X_p)]^2$ to be 
the normalized asymptotic variance from mutation for the strategy described above. 
We compare this to naive simulation 
($C_p^j \equiv 1$ for all 
$p$ and $j$ and $N_p \equiv N$ for all $p$) where by the law of large numbers,
\begin{equation*}
N{\mathbb E}\left[\left.(M_{p+1}-{\hat M}_{p})^2\right| {\mathcal F}_{p}\right] 
= \frac{1}{N}\sum_{j=1}^{N} 
 F(\xi_p^j)^2 \xrightarrow{a.s.} {\mathbb E}[F(X_p)^2]\quad \text{ as }N \to \infty,
\end{equation*}
so that ${\mathbb E}[F(X_p)^2]$ is the 
normalized asymptotic variance from mutation.
\end{remark}

\section{Binning}\label{sec:binning}

In traditional WE, the number 
of times $C_p^j$ we select 
particle $\xi_p^j$ is based 
on a binning technique. 
At each time step $p$, state space $E$ 
is divided 
into disjoint {\em bins} $B^r$, $r=1,\ldots,R$. 
That is, $E = \cup_{r=1}^R B^r$ where 
the union is disjoint. In general, 
the bins can be chosen adaptively; see 
Remark~\ref{rmk}. However, we 
will focus on fixed bins here and below.
In this setting, the selection step proceeds 
as follows. First, a 
target number of particles $N_p^r$ 
is set for each bin at time $p$. In many applications 
(see for instance~\cite{dan2,suarez2,suarez,dan1} and 
references), the target numbers are 
chosen so that the resulting 
particles cover space uniformly 
in some sense, which usually means $N_p^r \approx N/R$. We will take a different approach in the next section. The $C_p^j$'s are 
chosen such that either 
\begin{equation}\label{former}
\sum_{j:\xi_p^j \in B^r} {\mathbb E}[C_p^j|{\mathcal F}_p] = N_p^r,
\end{equation}
or, conditionally on ${\mathcal F}_p$, 
\begin{equation}\label{latter}
\sum_{j:\xi_p^j \in B^r} C_p^j = N_p^r.
\end{equation} 
In 
the latter case~\eqref{latter}, the number of particles in a given bin has a fixed deterministic value, while 
in the former~\eqref{former}, only the {average} number of 
particles in each bin is fixed. 
See Remark~\ref{rmk1} above. In~\eqref{former}, 
as discussed above, the $C_p^j$'s are usually chosen to have small variance, so the 
number of particles in a given bin has 
small variance. Throughout we will focus only on the case~\eqref{former}. Note 
that number of selected particles in bin $B^r$, namely 
$$\sum_{j:\xi_p^j \in B^r}C_p^j,$$ can 
equal zero. However, under assumption (A2), whenever there are particles 
in $B^r$ at time $p$, the 
expected number $N_p^r$ of selected particles in 
$B^r$ must be strictly positive. 
It is okay if there are no particles 
in a bin before selection -- that bin will just remain empty 
after the selection step.
\vskip5pt
\begin{remark}\label{rmk} We will use 
bins that are {fixed} in time. This 
stance allows us,  
in principle at least, to define a 
Markov state model~\cite{sarich,schutte} 
on the bins. This model can, 
in turn, be useful for minimizing 
the variance in~\eqref{eq_l2error}. 
We note, however, that the bins can be chosen adaptively 
and still fit the framework of Definition~\ref{def_PS}. 
For example, the bins can be deterministic functions 
of the particles and weights up to 
and including the current time. 
Theorem~\ref{thm_martingale} then demonstrates 
that WE samping is unbiased even when 
the bins are adaptively chosen. 
\end{remark}
\vskip5pt

\begin{algorithm}\caption{Constructing a coarse model}\label{alg3}
Choose bins $B^1,\ldots,B^R$ forming 
a partition of $E$, a sampling measure $\zeta$ 
on $(E,{\mathcal E})$, and a bounded 
function $f:E \to {\mathbb R}$. Then do the following.
\vskip5pt
\begin{itemize}
\item[1.] Estimate the probability ${P}_{rs}$
for $(X_p)_{p \ge 0}$ 
to go from $B^r$ to $B^s$ in one step:
\begin{equation*}
P_{rs} = \zeta(B^r)^{-1}\int_{B^r}{\mathbb P}\left[\left.X_{p+1} \in B^s\right|X_p = x\right]\,\zeta(dx).
\end{equation*}
Estimate the value of $f$ inside bin $B^r$ by $u_r$:
\begin{equation*}
u_r = \zeta(B^r)^{-1}\int_{B^r} f(x)\,\zeta(dx).
\end{equation*}
\item[2.] Let $v_p^r$ 
be the $r$th entry of the vector $P(P^{n-p-1}u)^2-(P^{n-p}u)^2$, where 
$u = (u_r)$ is considered a column vector and 
the squaring operations are entrywise. 
\vskip2pt
\item[3.] Let $\mu$ be the stationary 
distribution for the 
transition matrix $P = (P_{rs})$. 
That is, ${\mu} = (\mu_r)$ 
is the normalized left 
eigenvector of $P$ with eigenvalue $1$.
\end{itemize}
\vskip5pt
$P$ and $u$ in Step 1 can be obtained sampling many one-step trajectories 
starting at $\zeta$. A simple choice 
for 
$\zeta$ in a general setting would be the uniform (Lebesgue) 
measure. See also the Appendix for comments on another 
possibility for $\zeta$.
\end{algorithm}

\section{Adapting to a coarse model}\label{sec:adapt}

Suppose we have a coarse model for $(X_p)_{p 
\ge 0}$ and we want to use it to guide our sampling. 
As above, we fix $n\ge 0$ and a bounded function 
$f:E \to {\mathbb R}$. The coarse model 
will be adapted to some fixed choice of bins; 
we assume again that $E$ is divided into 
disjoint bins $B^r$, $r=1,\ldots,R$. 
We think of the coarse model as a Markov state model, where the states are the bins. More precisely, the coarse model will consist of  
approximations of the probability 
$P_{rs}$ that $X_{p+1} \in B^s$, given 
that $X_p \in B^r$, as well as estimates 
$u_r$ of the value of $f$ on $B^r$. 
Thinking of $P$ as a matrix and $u$ a column vector,  
let 
\begin{equation*}
v_p^r = r\text{th entry of the vector } P(P^{n-p-1}u)^2-(P^{n-p}u)^2,
\end{equation*}
where the squaring operations 
are entrywise. Then $v_p^r$ estimates the 
value in $B^r$ of
\begin{equation*}
Kg_{p+1}^2-g_p^2 \equiv K(K^{n-p-1}f)^2-(K^{n-p}f)^2.
\end{equation*}

Below we show how to use the coarse 
model to define a WE sampler 
so that $\eta_n(f)$ estimates ${\mathbb E}[f(X_n)]$ 
with small variance, using an approximate 
version of the strategy 
described in Section~\ref{sec:min}. Because we are 
using a coarse model that does not distinguish 
between 
points in a given bin, it 
is reasonable to 
take all the selected weights ${\hat \omega}_p^i$ 
of particles in a given
bin $B^r$ to have the same value ${\bar \omega}_p^r$: 
\begin{equation}\label{eq_reweight}
{\hat \omega}_p^i = {\bar \omega}_p^r, \qquad \text{if }{\hat \xi}_p^i  \in B^r.
\end{equation}
This type of weighting scheme 
is simply a choice of the practitioner. 
In particular, it leads to a class of WE samplers satisfying Definition~\ref{def_PS}, as follows.
In light of Definition~\ref{def_PS},  
the number of times $\xi_p^j \in B^r$ is selected is 
\begin{equation}\label{baromega}
{\mathbb E}[C_p^j|{\mathcal F}_p] = \frac{\omega_p^j}{{\bar \omega}_p^r}, \qquad \text{if }\xi_p^j \in B^r.
\end{equation}
Setting the average number of particles in $B^r$ as $N_p^r$ as in~\eqref{former}, we obtain
\begin{equation}\label{baromega2}
{\bar \omega}_p^r = \frac{\sum_{j:\xi_p^j \in B^r} \omega_p^j}{N_p^r}.
\end{equation}
Thus, the weighting scheme in~\eqref{eq_reweight}, together with 
a choice of target particle numbers $N_p^r$, 
leads to unique formulas 
for the selected weights 
and the expected number of children
of each particle.

From~\eqref{min_cpj}, the variance from 
mutation is minimized when 
\begin{align}\begin{split}\label{min_cpj2}
N_p^r &\equiv \sum_{j:\xi_p^j \in B^r} {\mathbb E}[C_p^j|{\mathcal F}_p] \\
&= \sum_{j:\xi_p^j \in B^r}\frac{N \omega_p^j \sqrt{Kg_{p+1}^2(\xi_p^j)-g_p(\xi_p^j)^2}}{\sum_{\ell=1}^{N_p}\omega_p^\ell \sqrt{Kg_{p+1}^2(\xi_p^\ell)-g_p(\xi_p^\ell)^2}}\\
&\approx \frac{N\sqrt{v_p^r}\sum_{j:\xi_p^j \in B^r}\omega_p^j}{\sum_{r=1}^R \sqrt{v_p^r} \sum_{j:\xi_p^j \in B^r}\omega_p^j}.
\end{split}
\end{align}
In Algorithms~\ref{alg2}-~\ref{alg4}, the $C_p^j$'s are independent with
\begin{equation}\label{Cpj}
C_p^{j} = \begin{cases} \lfloor \omega_p^j/{\bar \omega}_p^r\rfloor, & 
\text{w.p.}\quad 1- \left(\omega_p^j/{\bar \omega}_p^r-\lfloor \omega_p^j/{\bar \omega}_p^r\rfloor\right) \\ 
\lfloor \omega_p^j/{\bar \omega}_p^r\rfloor + 1, & \text{w.p.}\quad\omega_p^j/{\bar \omega}_p^r-\lfloor \omega_p^j/{\bar \omega}_p^r\rfloor\end{cases}, \qquad \text{if }\xi_p^j \in B^r,
\end{equation} 
where $\lfloor x \rfloor$ denotes the integer 
part of $x$, and the $N_p^r$'s are defined by
\begin{equation}\label{Npr2}
N_{p}^r := \frac{(N-{\tilde N}R)\sqrt{v_p^r}\sum_{i:\xi_p^i \in B^r}\omega_{p}^i}
{\sum_{r=1}^R \sqrt{v_p^r}\sum_{i:\xi_p^i \in B^r}\omega_{p}^i} + {\tilde N},
\end{equation}
where ${\tilde N}  \in (0,N/R)$ is a lower threshold 
for the target number of particles per bin, and by convention $N_p^r = {\tilde N}$ if the 
denominator in~\eqref{Npr2} is zero.
The $C_p^j$'s in~\eqref{Cpj} have been chosen to
minimize the variance due to selection over all possible 
choices satisfying~\eqref{baromega}.
See Algorithms~\ref{alg2} and~\ref{alg4} for implementations 
of these ideas. 
\begin{algorithm}\caption{A WE sampler adapted to a coarse model}\label{alg2}
Choose bins $B_r$, $r=1,\ldots,R$ forming 
a partition of $E$, a target total 
number of particles $N$, a lower threshold ${\tilde N}$, a final 
time $n$, and a bounded function $f:E \to {\mathbb R}$. Let $v_p^r$ be obtained as in Algorithm~\ref{alg3}. 
Choose initial points and weights with 
the distribution of $X_0$ in the sense of (A1). 
For $0 \le p \le\min\{n,\tau_{kill}\}$, iterate 
the following:
\vskip5pt
\begin{itemize}
\item[1.] Select $N_p^r$ according to~\eqref{Npr2} 
and define ${\bar \omega}_p^r$ as in~\eqref{baromega2}. 
\item[2.] Let $C_p^j$ be random variables defined 
by~\eqref{Cpj}, and select $\xi_p^j$ exactly $C_p^j$ times. 
Let ${\hat \xi}_p^i$, $i=1,\ldots,{\hat N}_p$ be the 
selected particles, with~${\hat N}_p = \sum_{j=1}^{N_p} C_p^j$.
\item[3.] Set $N_{p+1} = {\hat N}_p$ and 
assign the weight ${\omega}_{p+1}^i = {\bar \omega}_p^r$ if ${\hat \xi}_p^i \in B^r$.
\item[4.] Evolve the particles ${\hat \xi}_p^i$, $i=1,\ldots,{\hat N}_p$, independently
according to the law of $(X_p)_{p\ge 0}$ to get the next generation $\xi_{p+1}^i$,
$i=1,\ldots,N_{p+1}$ of particles.   
\end{itemize}
\vskip5pt
When $p=\min\{n,\tau_{kill}\}$, stop and output $\eta_n(f)$, an estimate of ${\mathbb E}[f(X_n)]$.
\end{algorithm}

Why did we set a lower threshold ${\tilde N}$ 
in~\eqref{Npr2}?
If ${\tilde N} = 0$ and $v_p^r$ is zero in a bin 
that contains particles, then
no particles can be selected in 
this bin and so
(A2) does not hold. 
Taking ${\tilde N}>0$ 
eliminates this problem 
by ensuring $N_p^r>0$ in each 
bin so that each particle has a 
positive survival probability.

Moreover, if $v_p^r$ is very small in some bins and large in others, if ${\tilde N} = 0$ then  
some selected particles can have very large weights due to a tiny value of $N_p^r$ 
in~\eqref{baromega2}.
While this is fine in principle -- the method 
is still unbiased --  we observed numerically 
that it is better to keep a target 
number of particles per bin that 
is bounded significantly away from zero. 
Note that, 
as desired, the expected number of selected 
particles is
$$\sum_{r=1}^R N_p^r = N,$$
unless $v_p^r = 0$ in every bin 
that contains particles, 
in which case we instead have $$ \sum_{r=1}^R N_p^r = {\tilde N}R.$$

If $\sum_{r=1}^R N_p^r < N_p$, then it is possible that ${\mathbb E}[C_p^j|{\mathcal F}_p] < 1$ for all 
$j$. If in addition the $C_p^j$'s 
are independent conditional 
on ${\mathcal F}_p$, then 
there is a strictly positive probability 
that all particles are killed 
in selection, 
that is, $\tau_{kill} = p+1$. 
However, we believe extinction is 
a remote possibility with an appropriate 
choice of parameters. Indeed in our simulations
we did not observe any events where all the particles were killed so 
long as $N$ was kept reasonably 
large and ${\tilde N}$ not too 
close to zero; see  Section~\ref{sec:example} below.

\section{Stationary averages}\label{sec:stationary}

In this section, let $(X_p)_{p \ge 0}$ have a 
unique stationary distribution $\pi$. Suppose 
we want to sample $\pi(f)$, where $f:E \to {\mathbb R}$ 
is bounded. Assume we have bins 
$B^r$, $r=1,\ldots,R$, and a coarse model for 
$(X_p)_{p \ge 0}$ as in the last section. This model 
consists of a coarse transition matrix $P$ which 
can be used to estimate variances as discussed above. 
Note that $P$ can also be used to estimate $\pi$. 
That is, the stationary vector $\mu$ of $P$ 
-- namely, the normalized left eigenvector for 
eigenvalue $1$ -- is a coarse estimate of 
$\pi$. 

To sample $\pi(f)$, we can begin with 
points approximately distributed according 
to $\pi$ in the some sense. Using the coarse model, and
beginning with $N$ points roughly uniformly distributed 
in space, we take initial points and weights with
\begin{equation}\label{eq_initial}
\#\{j:\xi_0^j \in B^r\} \approx \frac{N}{R}, \qquad \omega_0^j = \frac{\mu_r}{\#\{k:\xi_0^k \in B^r\}},\quad \text{if }\xi_0^j \in B^r.
\end{equation}
The final time $n$ should be large enough to 
allow the WE sampler to relax to the true 
stationary distribution $\pi$. Techniques which 
employ a coarse model to estimate $\pi$, 
and use this as an initial condition for WE, 
have appeared recently in~\cite{dan2,darve}. 
However, using the coarse model to minimize 
the variance in the above fashion appears to 
be new, to the best of our knowledge. 
In this context, minimizing the variance 
requires minimal additional work, since 
we already have the coarse transition matrix $P$ 
which can be used to estimate the quantities needed to minimize 
variance.

The question of how to choose the 
final time $n$ is difficult in general. 
Note, however, that we have some 
prior information from our coarse 
model. In particular, the 
second-largest (in absolute value) eigenvalue
$\lambda_2$ of $P$ can give us some idea of 
how fast convergence can be, from the 
heuristic
\begin{equation*}
{\mathbb E}[\eta_n(f)] - \pi(f) \approx O(\lambda_2^n).
\end{equation*}
Moreover the constant associated with the 
big O may be small due to the 
initial condition~\eqref{eq_initial}, 
though this is difficult to quantify 
without prior information about 
how close the initial condition is to $\pi$.
\vskip5pt
\begin{remark}
We comment briefly on two 
other possibilities for sampling $\pi(f)$. 

First, we could build the coarse model adaptively, using a Monte Carlo estimator of $v_p^r$.   
That is, we update $v_p^r$ at each 
time $p$ 
using the WE trajectory. 
One advantage of this is that, 
if we are using $\eta_n(f)$ 
to estimate $\pi(f)$, the most important 
contributions to the variance come 
from the final steps ($p$ near $n$), at which $v_p^r$ 
is the most accurate. 

Another possibility is, instead of 
estimating $\pi(f)$ from $\eta_n(f)$, 
we could use a time average via
$\pi(f) \approx (n+1)^{-1}(\eta_0(f)+\ldots+\eta_{n}(f))$. In this setting, 
it is natural to adaptively build 
estimates $v^r$ of $Kf^2-(Kf)^2$ in 
$B^r$, and plug them into~\eqref{Npr2} in place of the $v_p^r$ at each step. 
We do not test these strategies here, 
but leave them as interesting problems for future work.
\end{remark}

\section{Numerical example}\label{sec:example}

In this section $(Y_t)_{t \ge 0}$ is 
an Markov chain designed 
to mimic 
MD in a simple one dimensional 
energy landscape, and $X_p = Y_{p \delta t}$. 
In the context of WE, this means 
we resample from $(Y_t)_{t \ge 0}$ at each 
time interval $\delta t$. Our goal is 
to show that the adaptive 
sampling from the last section potentially
can be better than 
naive sampling or traditional 
WE sampling. Applying 
these ideas to more 
``realistic'' models in computational 
chemistry will be the 
focus of another work. 

\begin{algorithm}\caption{A WE sampler for stationary averages}\label{alg4}
Choose bins $B_r$, $r=1,\ldots,R$, a target total 
number of particles $N$, a final 
time $n$, and a bounded function $f:E \to {\mathbb R}$. Construct a coarse model as in Algorithm~\ref{alg3}. 
\vskip5pt
\begin{itemize}
\item[1.] Choose initial points and weights as in~\eqref{eq_initial}. 
\item[2.] For $0 \le p \le\min\{n,\tau_{kill}\}$, 
proceed through Algorithm~\ref{alg2}. 
\end{itemize}
\vskip5pt
The output $\eta_n(f)$ is an estimate of $\pi(f)$, the stationary average of $f$.
\end{algorithm}

More precisely, let  
$Y_t$ 
have values in $E = \{1,2,\ldots,90\}$ 
and transition matrix 
\begin{align*}
Q(i,i+1) &= \frac{2}{5}+\frac{m(i)}{5}, \qquad i=1,\ldots,89,\\
Q(i,i-1) &= \frac{2}{5}-\frac{m(i)}{5},\qquad i=2,\ldots,90,
\end{align*}
where 
\begin{equation*}
m(j) := \sin\left(\frac{6\pi j}{90}\right),
\end{equation*}
$Q(i,j) = 0$ if $|i-j|>1$ and $Q(i,i)$ is chosen 
so that $Q$ is stochastic. 
This is a discrete state space 
which mimics a one dimensional 
potential energy 
landscape with $3$ energy 
wells; see the bottom 
right of Figure~\ref{fig2}.
We take resampling intervals 
$\delta t = 4$, so
\begin{equation*}
X_p := Y_{4p}
\end{equation*}
and the transition matrix of $(X_p)_{p \ge 0}$ 
is $K = Q^4$. The resampling 
intervals are chosen so that 
a sufficiently large fraction 
of particles can change bins 
in each resampling time. 
The bins will be
\begin{equation*}
 B^r = \{3r-2,3r-1,3r\}, \qquad r = 1,\ldots,30.
\end{equation*}
Thus, there are $R = 30$ bins.
Let $\pi$ be the stationary 
distribution of $(X_p)_{p \ge 0}$, and 
\begin{equation*}
f(i) = \begin{cases}1, & 28 \le i \le 33 \\ 0, & \text{else}\end{cases}.
\end{equation*}
We also let ${\bar f} = f/6$ be 
its normalized version, which is 
useful for comparing with sampling 
distributions; see Figure~\ref{fig2} 
below.

We will 
obtain empirical approximations 
$\eta_n(f)$ of $\pi(f)$ for  
relaxation times~$n = 5,10,15,20,25,30$, 
using three 
types of sampling described below.
In all 
our simulations, we set a target 
number $N = 150$ of particles, 
Our initial 
particles and weights are 
the same for all simulations. 
They are chosen by constructing 
a coarse model as in Algorithm~\ref{alg3} 
with $\zeta$ the uniform 
measure on $E$, 
$\zeta(i) = 1/90$ for all $i \in E$.
Thus, our initial points and 
weights are chosen according to the distribution
\begin{equation*}
\nu_0(i) = \frac{\mu_r}{3},\qquad \text{if }i \in B^r.
\end{equation*}
In all our simulations we had $\tau_{kill}>n$. 

The first type of sampling 
uses  
Algorithm~\ref{alg4}, 
the procedure described above 
for adapting WE to a coarse model. We call this 
adaptive WE sampling. We 
construct a coarse model using Algorithm~\ref{alg3}
with the uniform sampling 
measure $\zeta$ described above. 
We use the parameters
$N = 150$ and ${\tilde N} = 1$.

In the second type of sampling 
we used a fixed 
target number of particles 
per bin. We 
call this traditional WE sampling. 
It is the sampling method described 
in~\cite{dan2,darve}. We 
use Algorithm~\ref{alg4} again, 
but instead of using a coarse model to define $N_p^r$ 
via~\eqref{Npr2}, we set $N_p^r = 5$ 
constant. This corresponds to distributing 
the $N = 150$ particles uniformly among the bins.

\begin{figure}\label{fig1}\vskip-190pt
  \centerline{\includegraphics[scale=.84]{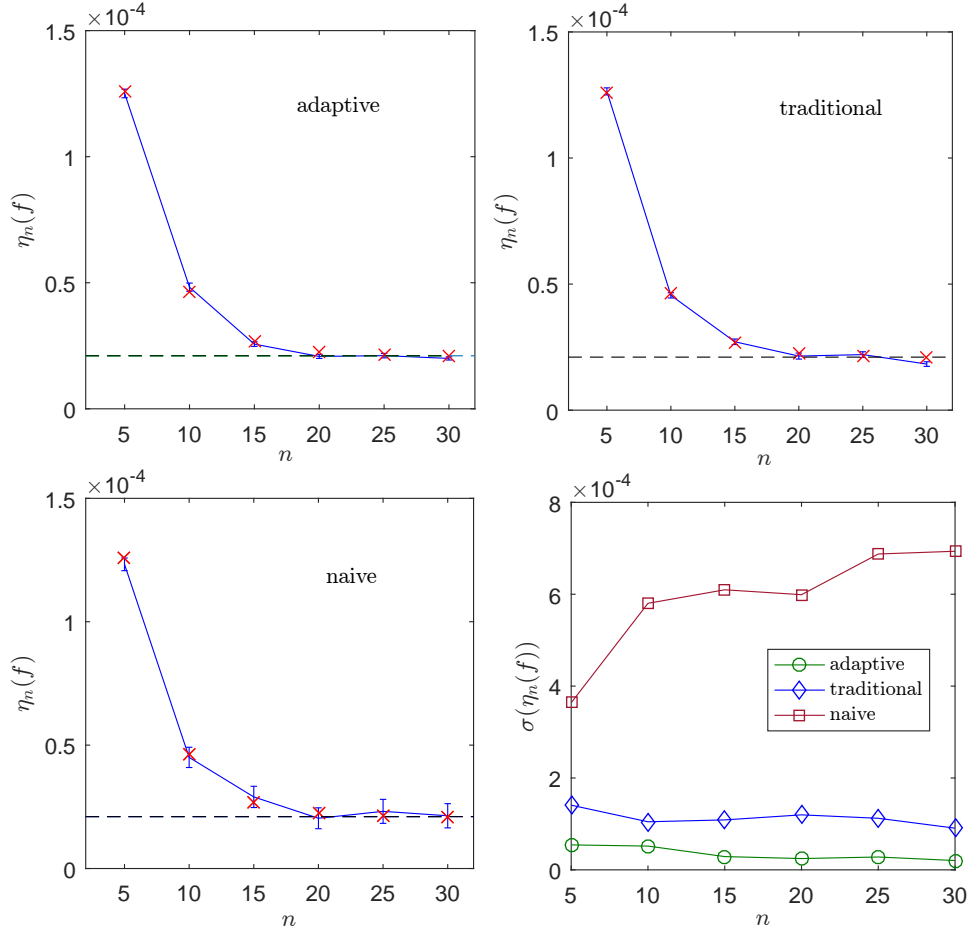}}
\vskip-150pt
\caption{Values of $\eta_n(f)$ 
vs. $n$ from the example 
in Section~\ref{sec:example} 
from adaptive, traditional 
and naive sampling. Data 
for adaptive, traditional 
and naive sampling is obtained 
from $10^3$, $10^4$, and $5 \times 10^4$ 
independent simulations, respectively. 
The crosses 
are exact values corresponding 
to $\nu_0 K^n f$, 
and the dotted line is 
the stationary value $\pi(f)$. {Bottom right}: 
sample standard deviations 
$\sigma(\eta_n(f))$ of 
$\eta_n(f)$, computed from the independent 
simulations. (The error bars in the 
other plots are these 
standard deviations divided by the 
square roots of the number of independent 
simulations.)}
\end{figure}

\begin{figure}\label{fig2}\vskip-195pt
  \centerline{\includegraphics[scale=.84]{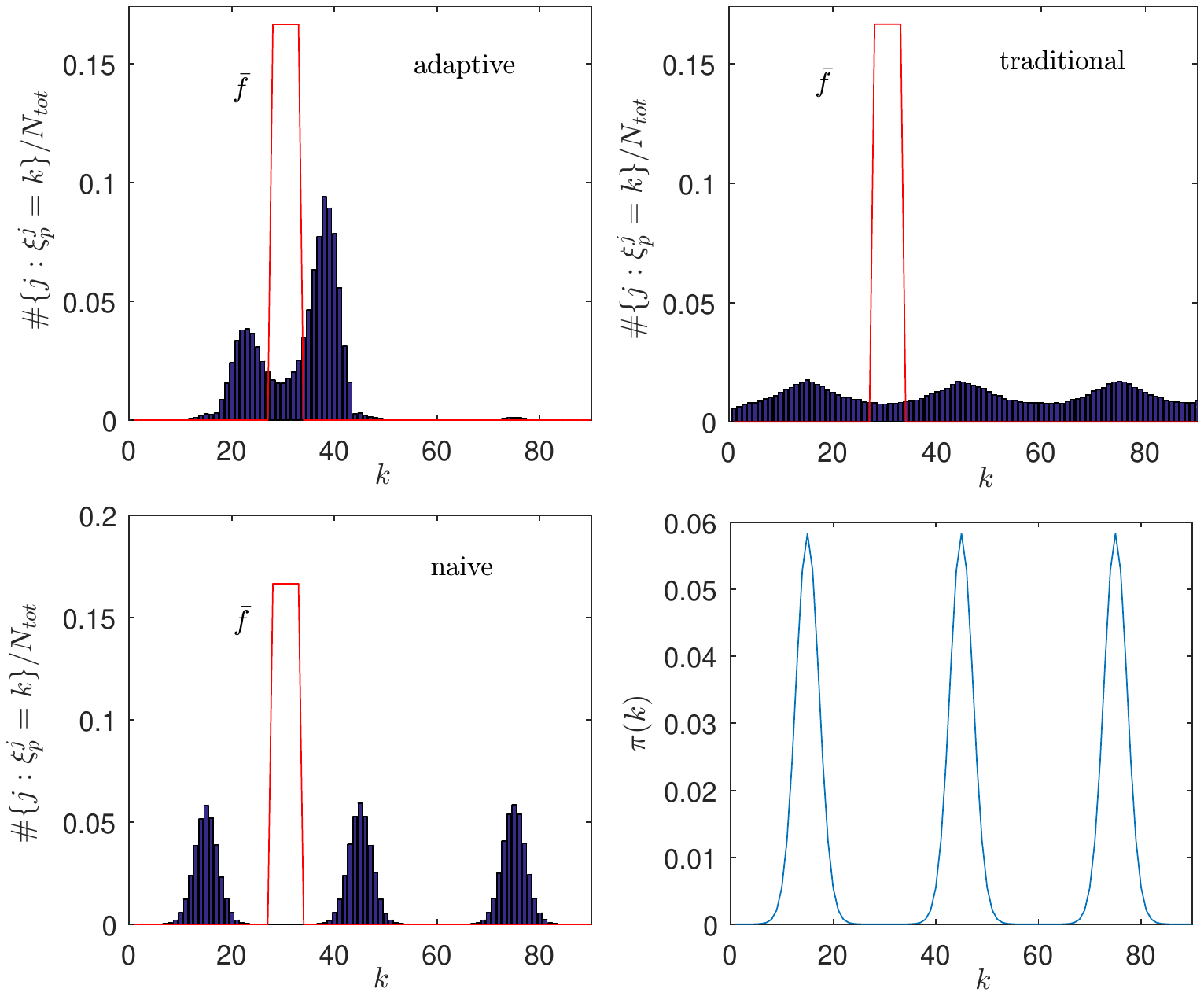}}
\vskip-170pt
\caption{Average distribution 
of particles $\xi_p^j$ at time $n = 30$, 
compared to ${\bar f}$, the normalized 
version of $f$. Here, $N_{tot}\equiv N_n$ is the 
total number of particles at time $n$.
{Bottom right}: 
Stationary distribution $\pi$.}
\end{figure}

The third type of sampling
does not use selection 
at all. We call this naive 
sampling. Here, we choose $N = 150$ initial particles and weights according to~\eqref{eq_initial}, and then
we simply evolve these particles 
independently until time $n$, 
without changing the weights. 

Results comparing adaptive WE sampling, 
traditional WE sampling and naive sampling are in 
Figures~\ref{fig1}-\ref{fig2}. In Figure~\ref{fig1}, we plot $\eta_n(f)$ 
vs. $n$ for various values 
of $n$, showing convergence 
to the stationary value 
$\pi(f)$. We compute error 
bars using
empirical standard deviations 
from $10^3$, $10^4$ and $5 \times 10^4$ independent 
simulations for adaptive, 
traditional, and naive sampling 
respectively. (We had to run more 
simulations for traditional WE
and naive sampling to get the 
numerics to converge.) The sample standard deviation for adaptive WE
sampling is significantly 
smaller than that of traditional WE
and naive sampling. 

In Figure~\ref{fig2}, 
we plot histograms representing 
the average distribution of 
the particles $\xi_n^i$ at time 
$n$. Note that traditional WE
sampling distributes the 
particles roughly uniformly in space,
as expected, while adaptive WE
sampling guides the particles 
towards the region 
in state space relevant 
for computing $f$.
Meanwhile, naive sampling 
distributes the particles 
approximately according to the stationary 
distribution $\pi$.

\begin{figure}\label{fig3}\vskip-210pt
  \centerline{\includegraphics[scale=.82]{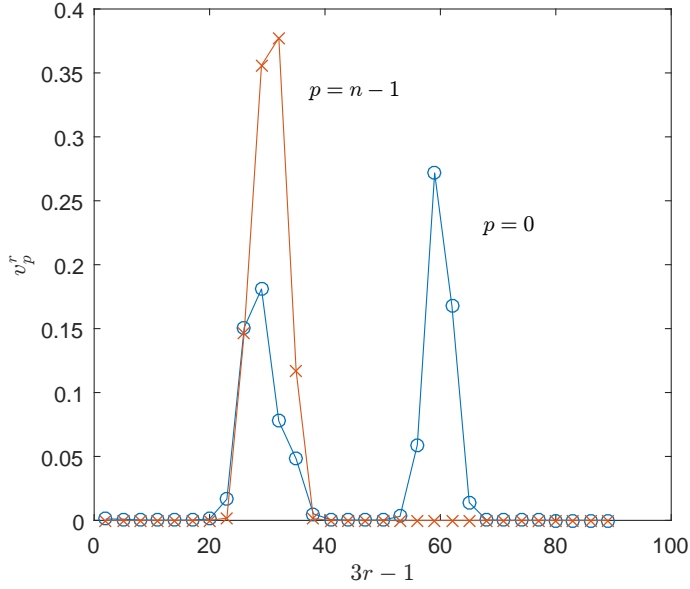}}
\vskip-210pt
\caption{The estimates $v_p^r$ vs. 
$p$ for the example 
in Section~\ref{sec:example}. 
Here, we take $p=0$ and $p = n-1$, 
where $n=30$ is the final time. 
}
\end{figure}

In Figure~\ref{fig3}, we 
plot the estimates 
$v_p^r$ from the adaptive 
sampling strategy for 
$p = 0$ and $p=n-1$ 
where $n=30$ is the relaxation time. Note that by time $n-1$, the sampling 
is focused near the support 
of $f$. 

When 
$f$ is a function with 
large values in regions of low $\pi$ 
probability, as in this 
example, naive sampling 
performs poorly compared 
to both traditional 
and adaptive WE sampling.
When state space is 
very large compared to 
the region $S$ where 
$f$ has large values 
(or is non-negligible), 
we expect adaptive WE
sampling to perform 
much better than traditional WE
sampling, due to the 
fact that traditional WE 
sampling will distribute 
the particles very thinly throughout 
space, including in $S$, while 
adaptive WE sampling will push 
most of the particles towards $S$.

A possible drawback of 
adaptive WE sampling is that it 
requires more computations 
at the resampling times, 
compared to traditional 
WE 
sampling. However, 
in practice the resampling 
times may be large enough 
so that this extra 
effort contributes 
little to the overall 
computational cost. 

Finally, we note that 
the adaptive sampling 
above can also be used 
more generally 
to estimate time 
marginals of $(X_p)_{p \ge 0}$, 
that is, expectations 
of the form ${\mathbb E}[f(X_n)]$ 
at fixed finite times $n$, 
from an arbitrary 
initial distribution of $X_0$. 
This is Algorithm~\ref{alg2}.
In this case, a MSM 
is still required to guide 
the sampling. One 
of the advantages of 
the adaptive sampling 
in the stationary case 
is that 
a MSM has already 
been computed as part 
of a preconditioning step.

\section*{Appendix: Computing dynamics from stationary averages}\label{sec:dynamics}

In this Appendix we show 
how to compute certain  
dynamical averages of $(X_p)_{p \ge 0}$ from 
stationary calculations. 
As above, $(X_p)_{p\ge 0}$ is a
time homogeneous Markov chain 
with values in~$E$.
The Hill relation~\cite{hill} 
shows that a mean hitting time 
can be reformulated as a 
certain stationary average. 
Similar ideas have recently 
been adapted to the time 
inhomogeneous  
setting; see~\cite{vankoten}. 
Here we focus on the time 
homogeneous case. 

By way of motivation, suppose 
we have a Markov chain with a 
(perhaps time reversible)
transition kernel $K_0$. Suppose 
we are interested in averages of 
the Markov chain, starting 
at a distribution $\rho$ and 
up to the hitting time $\tau_F$ 
of some set $F$ disjoint from 
the support of $\rho$. To 
compute such averages, we consider 
a modified, non-time reversible 
transition kernel $K$ constructed 
by setting $K = K_0$ outside $F$ 
and $K = \rho K_0$ inside $F$. 
Clearly, if we can sample from 
$K_0$, then we can also sample from $K$, 
simply by sampling from $K_0$ outside $F$ 
and then instantaneously restarting at $\rho$ each 
time we reach $F$. The following 
result recasts an average 
of the Markov chain with kernel 
$K_0$ starting at $\rho$ and up to time $\tau_F$ as a 
stationary average of the 
nonreversible Markov 
chain with kernel $K$.
\vskip5pt
\begin{theorem}\label{thm}
Suppose there is a set $F \subseteq E$ 
and a probability measure $\rho$ on $E$ with support 
disjoint from $F$ such that:
\vskip5pt
\begin{itemize}
\item[(B1)] The transition kernel $K$ of 
$(X_p)_{p \ge 0}$ satisfies ${\mathbbm{1}}_F(x)K(x,dy) = \rho K(dy)$,
\item[(B2)] With $\tau_F = \inf\{p>0\,:\,X_p \in F\}$, ${\mathbb E}^\rho[\tau_F]<\infty$ 
and ${\mathbb P}^{x}[\tau_F<\infty]=1$ $\forall x \in E$. 
\end{itemize}
\vskip5pt
Then for any bounded $g:E \to {\mathbb R}$, 
\begin{equation}\label{eq_hill}
{\mathbb E}^\rho\left[\sum_{p=1}^{\tau_F} g(X_p)\right] = 
\frac{\pi(g)}{\pi(F)},
\end{equation}
where $\pi$ is the unique stationary 
distribution of $(X_p)_{p\ge 0}$.
\end{theorem}
\vskip5pt
\begin{proof}
Assumptions (B1)-(B2) show that $(X_p)_{p \ge 0}$ has 
a unique stationary distribution $\pi$. Indeed, it can be checked (see~\cite{durrett}, Section~5.6) that 
\begin{equation*}
\pi(A) = \dfrac{{\mathbb E}^\rho\left[\sum_{p=1}^{\tau_F} {\mathbbm{1}}_A(X_p)\right]}{{\mathbb E}^\rho[\tau_F]}.
\end{equation*}
Thus, 
\begin{equation*}
\frac{\pi(g)}{\pi(F)} = \dfrac{{\mathbb E}^\rho\left[\sum_{p=1}^{\tau_F} g(X_p)\right]}{{\mathbb E}^\rho\left[\sum_{p=1}^{\tau_F} \mathbbm{1}_F(X_p)\right]} = {\mathbb E}^\rho\left[\sum_{p=1}^{\tau_F} g(X_p)\right].
\end{equation*}
\end{proof}
\vskip5pt

In practice, we are interested 
in the left hand side of~\eqref{eq_hill}. 
Assumption (B1) can be understood as  introducing 
a source $\rho$ and sink at $F$, 
while (B2) is an additional technical 
condition which ensures $\pi$ exists 
and is unique.  
In the context of the discussion 
above, (B1) corresponds to modifying 
the kernel $K_0$ of some underlying 
process to get the 
nonreversible kernel $K$. {This 
modification is 
only a computational 
tool}, as it does not affect the
LHS of~\eqref{eq_hill}. 

Thus, though the process we 
are interested in usually
does not satisfy (B1), 
we can modify it in $F$ 
so that (B1) holds, and 
meanwhile the left hand side of~\eqref{eq_hill} 
is the same for both the 
original and modified process. 
In this 
setting, if the original process is reversible, it 
is natural to take the sampling measure $\zeta$ in Algorithm~\ref{alg3} to 
be its stationary distribution, provided it can be efficiently 
calculated by Markov chain Monte Carlo or other common sampling techniques for 
reversible processes. 
It is important to note that 
such techniques cannot be used
to directly sample $\pi$, 
since the 
modified process is nonreversible.

Two special cases of~\eqref{eq_hill} are of particular 
interest. First, suppose $F = A \cup B$ is a disjoint 
union, $g = {\mathbbm{1}}_B$, and $\tau_S = \inf\{p>0:X_p \in S\}$ is the 
first time to hit $S$. Then
\begin{equation}\label{hit_prob}
{\mathbb P}^\rho\left[\tau_B < \tau_A \right] = \frac{\pi(B)}{\pi(A \cup B)}.
\end{equation}
Next, suppose $g \equiv {\mathbbm{1}}$. Then 
\begin{equation}\label{MFPT}
{\mathbb E}^\rho[\tau_F] = \frac{1}{\pi(F)}.
\end{equation}
Equation~\eqref{MFPT} is known as the Hill relation~\cite{hill}.
Equations~\eqref{hit_prob} and~\eqref{MFPT} 
show how stationary 
calculations can be used to compute hitting probabilities and hitting 
times. We can compute the right 
hand side of~\eqref{hit_prob} 
by applying Algorithm~\ref{alg4} above 
to $f = {\mathbbm{1}}_B$ and then 
$f = {\mathbbm{1}}_{A \cup B}$. 
Similarly, we can compute the right hand side of~\eqref{MFPT}
by applying
Algorithm~\ref{alg4} with $f = {\mathbbm{1}}_F$. A simple choice for $\rho$ would be 
$\rho = \delta_x$, the delta distribution 
at a point $x \notin F$. 
A more complicated but important case  
is the so-called  
equilibrium hitting time 
between an initial set $I$ and final set $F$; see for instance~\cite{Bhatt} for definitions and discussion. 
In this case, $\rho$ is the 
distribution of endpoints of trajectories 
under the original kernel $K_0$ stopped upon hitting $I$ and which last came from $F$. Sampling this 
distribution can be difficult in general~\cite{Bhatt}.

We conclude by briefly connecting the 
discussion above to Exact Milestoning~\cite{aristoff,elber}, an 
algorithm mentioned in the 
Introduction for 
sampling dynamical quantities 
like mean hitting times. Consider 
the following seemingly more 
general framework. Suppose that $(Y_t)_{t \ge 0}$ 
is some underlying process and $(\tau_p)_{p \ge 0}$ are increasing stopping times for $(Y_t)_{t \ge 0}$ such that $(X_p)_{p \ge 1}$ 
defined by
$$X_p = (Y_{\tau_{p-1}+1},\ldots,Y_{\tau_{p}})$$ 
is a time homogeneous Markov chain in 
$\cup_{n=1}^\infty E^n$. For instance, if 
$(Y_t)_{t \ge 0}$ is a 
time homogeneous Markov chain,
we could take 
$\tau_p = p \delta t$ with $\delta t$ a 
deterministic time, as in the example 
in Section~\ref{sec:example}, or $\tau_p = \inf\{t > \tau_{p-1}: Y_t \in S\}$ for some set $S\subseteq E$, and $\tau_0 = 0$. 
The latter choice corresponds to Exact 
Milestoning, in 
which $S$ corresponds to the union 
of all the milestones. In this setting, if we take $g(X_p) = \tau_p-\tau_{p-1}$, 
$F = \cup_{n=1}^\infty (E^{n-1} \times R)$,
and $$T_R = \inf\{\tau_p > 0: Y_{\tau_p} \in R\},$$ then from~\eqref{eq_hill},
\begin{equation}\label{MFPT2}
{\mathbb E}^\rho[T_R] 
= \frac{{\mathbb E}^\pi[\tau_1]}{\pi(F)}.
\end{equation}
This is the equation on which Exact Milestoning is based; see for instance Theorem~3.4 of~\cite{aristoff}. Thus in Exact Milestoning, 
we can find the time $T_R$ for 
$Y_t$ to first reach $R$ starting at $\rho$ by 
computing $\pi(F)$ along with short 
trajectories of $Y_t$ starting 
at $\pi$ up to the first 
time to hit $S$.

\section*{Acknowledgments}
D. Aristoff would like to acknowledge 
enlightening conversations 
with Tony Leli\`evre, Petr Plech\'{a}\v{c}, Mathias 
Rousset, Gideon Simpson, 
Ting Wang, and Dan Zuckerman. 
The author also wishes to thank two 
anonymous referees for providing 
numerous useful suggestions for 
improving the manuscript. 
D. Aristoff also gratefully 
acknowledges support from the 
National Science Foundation 
via the award
NSF-DMS-1522398.

\end{document}